\documentclass[12pt]{article}
\usepackage{amsfonts,amsmath,amssymb}
\usepackage{enumerate}
\usepackage{theorem}
\numberwithin{equation}{section}

\newtheorem{theorem}{Theorem}[section]

\newtheorem{definition}{Definition}[section]

\newcommand{\cl}[1]{\mathcal{#1}} 
\newcommand{\bb}[1]{\mathbb{#1}}

\begin{document}

\newcommand{\ma}[1]{\marginpar{\tiny #1}}

\title{Some notes on Morita equivalence of operator algebras}

\author{G.K. Eleftherakis}
\date{}
\maketitle

\begin{abstract} In this paper we present some key moments in the history of 
Morita equivalence for operator algebras.
\end{abstract}

\bigskip



\bigskip
\bigskip

We briefly describe some events in the history of Morita equivalence for operator algebras.
We make no claim of completeness in this presentations. 
Also, we have modified some of the definitions so as to be compatible with our presentation.
The reader can use the books \cite{bm, paul} 
for the notions of operator space theory used 
in this article. In what follows if $H, K$ are Hilbert spaces, 
we denote by $B(H, K)$ the space of bounded operators from $H$ to $K$ and we write $B(H)$ for $B(H,H)$.

\section{ Morita equivalence of $C\sp{*} $-algebras}

M. A. Rieffel introduced the idea of Morita equivalence of $C\sp{*} $ and  $W\sp{*} $ algebras. 
 The definitions and results of 
this section can be deduced from \cite{bgr, rief, rief1, rief3}. Let $H, K$ be Hilbert spaces and $\cl T$ be a subspace of $B(H, K).$ The 
space $\cl T$ is called a ternary ring of operators (TRO) if $$\cl T\cl T^*\cl T\subset \cl T.$$ 

\begin{definition}\label{11} Let $A$ be a $C\sp{*} $ algebra. We denote by $ \;_A \mathfrak{HM} $ the category of representations of $A:$ 
The objects of $\;_A \mathfrak{HM}$ are pairs $(\pi , H)$ where $H$ is a Hilbert space and 
$$ \pi : A\rightarrow B(H) $$ is a $*-$homomorphism such that $\overline{\pi (A)(H)}=H.$ If $(\pi _i, H_i)\in \:_A \mathfrak{HM} , i=1,2$ 
the corresponding space of category homorphisms 
 is the following:
$$B_A(H_1, H_2)=\{x\in B(H_1, H_2): x\pi_1(a) =\pi _2(a)x \;\;\forall \;\;a\in A\}.$$
\end{definition} 

\begin{definition}\label{12} Let $A, B$ be $C\sp{*} $ algebras. These algebras are called strongly Morita equvalent in the sense of Rieffel (R-SME), if there exist 1-1  $*-$homomorphisms 
$$\pi : A\rightarrow B(H),\quad \rho  : B\rightarrow B(K) $$ 
where 
$H, K$ are Hilbert spaces and there exists a TRO $\cl T\subset B(H, K)$ such that 
$$\pi (A)=[\cl T^*\cl T] ^{-\|\cdot\|} ,\quad \rho(B) =[\cl T\cl T^*] ^{-\|\cdot\|} .$$
\end{definition}

\begin{theorem}\label{13} Let $A, B$ be R-SME  $C\sp{*} $ algebras. 
Then the categories $\;_A \mathfrak{HM}$ and $ \;_B \mathfrak{HM}$ 
are equivalent. 
\end{theorem}

\begin{definition} Let $\cl K$ be the space of compact operators acting on an infinite dimensional separable Hilbert space. The 
$C\sp{*} $ algebras $A$ and $B$ are called strongly stably isomorphic
 if the $C\sp{*} $ algebras $\cl A\otimes \cl K$ and $\cl B\otimes \cl K$  
are $C\sp{*} $  isomorphic. (The tensor product $\otimes $ is the spatial tensor product, 
see  \cite{bm}).
\end{definition}

We can easily see that if $A$ and $B$ are strongly, stably isomorphic then they are R-SME. The converse doesn't hold. But if 
$A$ and $B$ possess countable approximate identities the converse is true. In particular, if $A$ and $B$ are separable or unital, 
then they are R-SME if and only if they are strongly stably isomorphic.

\section{Morita equivalence of $W\sp{*} $-algebras}

\begin{definition}\label{21} Let $M$ be a $W\sp{*} $ algebra. We denote by $ \;_M \mathfrak{NHM} $ the category of normal representations 
of $M:$ 
The objects of $\;_M \mathfrak{NHM}$ are pairs $(\pi , H)$ where $H$ is a Hilbert space and 
$$ \pi : M\rightarrow B(H) $$ is a weak* continuous $*-$homomorphism such that the 
identity operator belongs to  $\pi (M).$ 
If $(\pi _i, H_i)\in \;_M \mathfrak{NHM} , i=1,2$ 
the corresponding space of category homorphisms is the following:
$$ B_M(H_1, H_2) =\{x\in B(H_1, H_2): x\pi_1(a) =\pi _2(a)x \;\;\forall \;\;a\in M\}.$$
\end{definition}

\begin{definition}\label{22} Let $M, N$ be  $W\sp{*} $ algebras. These algebras are called weakly Morita equvalent in the sense of Rieffel (R-WME), 
if there exist 1-1  weak* continuous $*-$homomorphisms 
$$\pi : M\rightarrow B(H), \quad \rho  : N\rightarrow B(K) $$ where 
$H, K$ are Hilbert spaces and there exists a TRO $\cl T\subset B(H, K)$ such that 
$$\pi (M)=[\cl T^*\cl T] ^{-w^*}, \quad \rho(N) =[\cl T\cl T^*] ^{-w^*} .$$
\end{definition}

\begin{definition}\label{23} Let $M, N$ be a $W\sp{*} $ algebras and $$\cl F: \;_M \mathfrak{NHM} \rightarrow \;_N \mathfrak{NHM} $$  
be a functor. This functor is called a $*$-functor if for every pair of 
objects $H_1, H_2\in \;_M \mathfrak{NHM} $ and every 
$x\in  B_M(H_1, H_2)$, the map $ \cl F(x^*) $  is equal to $\cl F(x)^*.$
\end{definition}

\begin{theorem}\label{24} Let $M, N$ be $W\sp{*} $ algebras. The following are equivalent:

(i) The algebras $M$ and $N$ are R-WME.

(ii) The categories $\;_M \mathfrak{NHM} , \;_N \mathfrak{NHM}$ are equivalent through a $*-$functor.

(iii) There exists a cardinal $I$ such that the algebras $M_I(M), M_I(N)$ of all 
$I\times I$ matrices whose finite 
submatrices have uniformly bounded norms are isomorphic as $W\sp{*} $ algebras.

(iv)  There exist faithful normal representations $\pi $ of $M$ and $\rho $ of $N$ such that the commutants $\pi(M)^\prime$ and  $ \rho (N)^\prime$ are isomorphic.
\end{theorem}

The equivalence of (i) and (ii) was proved in \cite{rief2}. The others can be found in \cite{bm}.

\section{Strong Morita equivalence of operator algebras.}

In this section we assume that all operator algebras (see \cite{bm, paul}) 
have contractive approximate identities (cai). 
The definitions and results of 
this section can be deduced from \cite{blecher, bmp}.

\begin{definition} \label{31} Let $A$ be an operator algebra and $X$ be an operator space. We say that $X$ is a left $A$-operator module 
if there exists a completely contractive bilinear map $$A\times X\rightarrow X. $$ Similarly we can define right $A$-operator modules.
\end{definition}

\begin{definition}\label{32} Let $A$ be an operator algebra. We denote by $ \;_A \mathfrak{OM} $ the category of left operator modules: 
The objects of  $\;_A \mathfrak{OM}$ are the left operator modules $X$ for which $$e_ix\rightarrow x\;\;\forall \;\;x\in X, $$ 
where $(e_i)_i$ is a cai for the algebra $A.$ 
If $X_1, X_2\in \;_A \mathfrak{OM}$, 
the corresponding space of  category homorphisms 
is the following: $$CB_A(X_1, X_2)=\{u: X_1\rightarrow X_2: u\;\;\mbox{is a compl. bounded $A$-module map } \}.$$
\end{definition}

\begin{definition}\label{33} Let $A, B$ be operator algebras. A functor $$\cl F: \;_A \mathfrak{OM} \rightarrow \;_B \mathfrak{OM} $$ is 
called completely contractive if the map $x\rightarrow \cl F(x)$ is completely contractive for all pairs 
$X_1, X_2\in \;_A \mathfrak{OM}$ and all $x\in CB_A(X_1, X_2)$. 
\end{definition}
 
\begin{definition}\label{34} Let $A, B$ be operator algebras. These algebras are called strongly Morita equivalent, in the sense of 
Blecher Muhly and Paulsen, (BMP-SME) if there exist completely isometric homomorphisms $$\pi : A\rightarrow B(H), \;\;\rho : 
B\rightarrow B(K),$$ where $H, K$ are Hilbert spaces such that $$\overline{\pi (A)(H)}=H, \;\;\;\overline{\rho (B)(K)}=K$$ 
and such that there exist a $(\rho(B),\pi(A))$-module $U$ and a 
$(\pi(A),\rho (B))$-module $V$ satisfying 
$$\pi (A)=[VU] ^{-\|\cdot \|} , \;\;\;\rho(B) =[UV]^{-\|\cdot \|} .$$ We also require 
the existence of nets of positive integers $(n_t)_t, (m_i)_i$ and nets 
$$(v_t^1)\subset Ball(R_{n_t}(V)), \;\;\;(u_t^1)\subset Ball(C_{n_t}(U))$$
$$(u_i^2)\subset Ball(R_{m_i}(U)), \;\;\;(v_i^2)\subset Ball(C_{m_i}(V))$$
such that $(v_t^1u_t^1)_t$ is an approximate identity for  $\pi(A) $ and $(u^2_iv_i^2)_i$ is an approximate identity for $\rho(B). $
\end{definition}

We recall (see \cite{bm, paul}) the Haagerup tensor product $U\otimes ^hV$ of
two operator spaces $U$ and $V.$ This space is an operator space and it 
has the property of linearizing completely bounded bilinear maps. Suppose $U$ is a right $A$ operator module and $V$ is a left $A$ operator module 
over an operator algebra $A.$  We denote by $\Omega $ the space 
$$[u\otimes (av)-(ua)\otimes v: \;\; u\in U\;\;v\in V\;\;a\in A ]^{-\|\cdot\|}$$ and 
by $U\otimes ^h_AV$ the quotient $(U\otimes ^hV)/\Omega$. 
This space has the property of linearizing $A$-balanced 
completely bounded bilinear maps. This means that if $X$ is an operator space and $$\phi : U\times V\rightarrow X$$ is a 
completely bounded bilinear map satisfying 
$$\phi(ua, v) =\phi (u, av)\;\;\forall \;\;u\in U,\;\;v\in V,\;\;a\in A $$ 
then  there exists a completely bounded linear map 
$$\hat \phi : U\otimes ^h_A V\rightarrow X$$ 
such that 
$$\hat \phi (u\otimes _Av)=\phi (u,v).$$
 If $B$ is an operator algebra and $U$ (resp. $V$) is a left (resp. right)  $B$ 
operator module then $U\otimes ^h_AV$ is a left (resp. right)  $B$ 
operator module. 

\begin{theorem}\label{35} Let $A, B$ be operator algebras. The following are equivalent:

(i) The algebras $A$ and $B$ are BMP-SME.

(ii) The categories $\;_A \mathfrak{OM},  \;_B \mathfrak{OM}$ are equivalent through a completely contractive functor. 

(iii) There exist a $(B,A)$ operator module $U$ and an  $(A,B)$ operator module $V$ such that $U\otimes ^h_AV$ and $B$ (resp.
 $V\otimes ^h_BU$ and $A$) are equivalent as $B$ (resp. $A$) operator modules. 
\end{theorem}

\begin{theorem}\label{36} Two $C\sp{*} $ algebras are R-SME if and only if 
they are BMP-SME.
\end{theorem}

\section{Weak* Morita equivalence of dual operator algebras.}

In this section we assume that all dual operator algebras (see \cite{bm}) are unital. 
If $M$ is a dual operator algebra, 
a normal representation of $M$  is a map $\alpha : M\rightarrow B(H),$ where $H$ is a Hilbert space and $\alpha $ is 
a completely contractive weak* continuous homomorphism.  The definitions and results of 
this section can be deduced from \cite{bk, kash}. 

\begin{definition} \label{41} Let $M$ be a dual operator algebra and $X$ be a dual operator space. We say that $X$ is a left dual $M$-operator 
module if there exists a completely contractive weak* continuous bilinear map $$M\times X\rightarrow X.$$ Similarly we can define  
right dual operator modules over $M.$
\end{definition}

\begin{definition}\label{42} Let $M$ be a dual operator algebra. We denote by $\;_M\mathfrak{NOM}$ the category of 
left dual operator modules: The objects of  $\;_M\mathfrak{NOM}$ are the left dual operator modules $X$ for which 
$\overline{MX}^{w^*}=X.$ If $X_1, X_2\in \;_M\mathfrak{NOM}$ the corresponding space $CB_M^\sigma (X_1, X_2)$ of  category 
morphisms is the space of all weak* continuous completely bounded $M$-module maps from $X_1$ into $X_2.$
\end{definition}

\begin{definition}\label{43} Let $M, N$ be dual operator algebras. These algebras are called weakly* Morita equivalent in the sence 
of Blecher and Kashyap (BK-WME) if there exist completely isometric normal unital  representations $$\pi: M\rightarrow B(H), 
\rho : N\rightarrow  B(K), $$ where $H, K$ are Hilbert spaces 
such that there exist a $(\rho (N),\pi( M))$ bimodule $U\subset B(H, K)$ and a 
$(\pi( M),\rho (N))$ bimodule $V\subset B(K,H)$ satisfying 
$$\pi (M)=[VU] ^{-w^*},\;\;\;\rho (N)=[UV] ^{-w^*}.$$ 
Also we require the existence of nets of positive integers $n_t, m_i$ and nets 
$$(v_t^1)\subset Ball(R_{n_t}(V)), \;\;\;(u_t^1)\subset Ball(C_{n_t}(U))$$
$$(u_i^2)\subset Ball(R_{m_i}(U)), \;\;\;(v_i^2)\subset Ball(C_{m_i}(V))$$
such that the nets $(v_t^1u_t^1)_t$ and $(u^2_iv_i^2)_i$ converge in the weak* topology to the identity operators. 
\end{definition}

We recall (see \cite{bm}) the normal Haagerup tensor product $U\otimes ^{\sigma h}V$ of the dual operator spaces $U$ and $V.$ 
This space is a dual  operator space and it 
has the property of linearizing weak* continuous completely bounded bilinear maps. 
Suppose $U$ is a right $M$ dual operator module and $V$ is a left $M$ dual operator module 
over a dual  operator algebra $M.$  We denote by $\Omega^\sigma  $ the space 
$$[u\otimes (av)-(ua)\otimes v: \;\; u\in U\;\;v\in V\;\;a\in A ]^{-w^*}$$ and 
by $U\otimes ^{\sigma h}_MV$ the quotient $(U\otimes ^{\sigma h}V)/\Omega^\sigma .$
 This space has the property of linearizing $M$-balanced weak* continuous 
completely bounded bilinear maps \cite{elepaul}. 
This means that if $X$ is a dual  operator space and 
$$\phi : U\times V\rightarrow X$$ 
is a weak* continuous 
completely bounded bilinear map satisfying 
$$\phi(ua, v) =\phi (u, av)\;\;\text{for all }\;\;u\in U,\, v\in V\;\;\text{and }\;a\in A $$
 then 
there exists a weak* continuous completely bounded linear map 
$$\hat \phi : U\otimes ^{\sigma h}_M V\rightarrow X$$ 
such that 
$$\hat \phi (u\otimes _Mv)=\phi (u,v).$$ 
If $N$ is an operator algebra and the above $U, $ (resp. $V$) is a left (resp. right )  
$N$ operator module then $U\otimes ^{\sigma h}_MV$ is a left (resp. right )  $N$ 
operator module \cite{elepaul}. 

\begin{definition}\label{43} Let $N, M$ be dual operator algebras. A functor $$\cl F: \;_M\mathfrak{NOM}\rightarrow \;_N\mathfrak{NOM} $$ 
is called completely contractive and normal if the map $$\cl F: CB_M^\sigma (X_1, X_2) \rightarrow CB_N^\sigma (\cl F(X_1), \cl F(X_2)) $$ 
is completely contractive and weak* continuous for all $X_1, X_2\in \;_M\mathfrak{NOM}.$
\end{definition}

\begin{theorem}\label{44}Let $M, N$ be dual operator algebras. The following are equivalent:

(i) The algebras $M$ and $N$ are BK-WME.

(ii) The categories $\;_M \mathfrak{NOM},  \;_N \mathfrak{NOM}$ are equivalent through a completely contractive and normal functor. 

(iii) There exist a $(N,M)$ dual operator module $U$ and an  $(M,N)$ dual operator module $V$ such that $U\otimes ^{\sigma h}_MV$ and $N$ (resp.
 $V\otimes ^{\sigma h}_NU$ and $M$) are equivalent as $N$ (resp. $M$) dual operator bimodules. 

\end{theorem}

\begin{theorem} \label{45} Two $W\sp{*} $ algebras are R-WME if
and only if they are BK-WME.
\end{theorem}

\begin{theorem} \label{46} Let $A, B$ be BMP-SME approximately unital operator algebras. Then their second duals $A^{**}, B^{**}$ 
are BK-WME algebras.
\end{theorem}

\begin{theorem} \label{47} \cite{elenest} Let $\cl N_1, \cl N_2$ be nests  (see \cite{dav}) 
acting on the Hilbert spaces $H_1, H_2$ respectively, let $M_1, M_2$ be 
the corresponding nest algebras and let $\cl K(M_1), \cl K(M_2)$  be the subalgebras 
of compact operators. The following are equivalent:

(i) The algebras $M_1, M_2$ are BK-WME.

(ii) The algebras $\cl K(M_1), \cl K(M_2)$ are BMP-SME.

(iii)The nests $\cl N_1, \cl N_2$ are isomorphic.

(iv) There exist an $(M_1, M_2)$ bimodule $U\subset B(H_2, H_1)$ and 
an $(M_2, M_1)$ bimodule $V\subset B(H_1, H_2)$ 
such that $$ M_1=[UV]^{-w^*}, M_2=[VU]^{-w^*}. $$ 
\end{theorem}

\section{Stable isomorphism of dual operator algebras.}

\begin{definition} \label{51}\cite{eletro} Let $M, N$ be weak* closed algebras acting on the Hilbert spaces $H$ and $K$ respectively. We 
call them TRO equivalent if there exists a TRO $\cl T\subset B(H, K)$ such that 
$$ M=[\cl T^*N\cl T]^{-w^*} \;\;\;\;N=[\cl TM\cl T^*]^{-w^*} .$$
\end{definition}

\begin{theorem}\label{52}\cite{eletro} 
Let $\cl L_1, \cl L_2$ be reflexive lattices and let $M=\mathrm{Alg}(\cl L_1)$, $N=\mathrm{Alg}(\cl L_2)$ 
 be the corresponding algebras of operators leaving invariant every element of 
$\cl L_1$ and $ \cl L_2$ (see \cite{dav}); also let 
$\Delta (M) = M\cap M^*$ and $\Delta (N)=N\cap N^*$ be their diagonals. 
The algebras $M$ and $N$ are TRO equivalent if and only if 
there exists a $*$-isomorphism 
$$\theta : \Delta (M) ^\prime=\cl L_1^{\prime \prime} 
\rightarrow \Delta (N) ^\prime=\cl L_2^{\prime \prime} $$ such that 
$\theta (\cl L_1)=\cl L_2.$
\end{theorem}

\begin{definition}\label{elepure} Let $M, N$ be dual operator algebras. We call them $\Delta$-equivalent if there exist completely isometric normal 
representations $\alpha $ of $M$ and $\beta $ of $N$ such that $\alpha (M)$ and $\beta (N)$ are TRO equivalent.
\end{definition}

As in definition \ref{21} if $M$ is a unital dual operator algebra, we denote by $\;_M\mathfrak{NHM}$ the category of normal 
completetely contractive  representations 
of $M$: If $(H_i, \alpha _i), i=1,2$ are objects of $\;_M\mathfrak{NHM}$, 
the space of morphisms  $Hom_M(H_1, H_2) $ is the following:
$$ Hom_M(H_1, H_2) =\{x\in B(H_1, H_2): x\alpha_1(m) =\alpha_2(m)x\;\;\forall \;\;m\in M \}.$$
Let $\Delta (M)=M\cap M^* $ be the diagonal of $M.$ Observe that $\alpha _i|_{\Delta (M)}$ is a $*-$homomorphism since $\alpha _i$ 
is a contraction. 

We also define the category $\;_M\mathfrak{DNHM}$ which has the same objects as $\;_M\mathfrak{NHM}$ but for every pair of objects $(H_i, \alpha _i), i=1,2$ 
 the space of morphisms is the following: 
$$ Hom_M^D(H_1, H_2) =\{x\in B(H_1, H_2): x\alpha_1(m) =\alpha_2(m)x\;\;\forall \;\;m\in 
\Delta (M) \}.$$ Observe that $$Hom_M(H_1, H_2) \subset Hom_M^D(H_1, H_2) .$$
We say that $\cl F: \;_M\mathfrak{DNHM}\rightarrow \;_M\mathfrak{DNHM} $ is a 
$*$-functor if for every pair $H_1, H_2\in\;_M\mathfrak{DNHM}$  
$$\cl F(x^*)=\cl F(x)^*\;\;\forall \;\;x\in  Hom_M^D(H_1, H_2) .$$

The main theorem of this section, which is a generalization of Theorem \ref{24},
 is the following:

\begin{theorem}\label{55}\cite{elepure, eletro, elepaul} Let $M, N$ be unital dual operator algebras. The following are equivalent:

(i) The algebras $M$ and $N$ are $\Delta -$equivalent .

(ii) There exists a $*$-functor of equivalence $$\cl F: \;_M\mathfrak{DNHM}\rightarrow  \;_N\mathfrak{DNHM} $$ such that 
$$\cl F(\;_M\mathfrak{NHM} )=\;_M\mathfrak{NHM} .$$

(iii) The algebras $M$ and $N$ are stably isomorphic in the following sense: 
there exists a cardinal $I$ such that 
 the algebras $M_I(M), M_I(N)$ of  $I\times I$ matrices whose finite 
submatrices have uniformly bounded norms are isomorphic as dual operator algebras.

(iv) There exist completely isometric normal representations $\alpha $ of $M$ and $\beta $ of $N$ such that 
$$\alpha(M)=  \mathrm{Alg}(\cl L_1) , \;\;\;\beta (N)=\mathrm{Alg}(\cl L_2) $$ 
for reflexive lattices $\cl L_1$ and $\cl L_2$ and there exists a $*-$isomorphism $$\theta : \alpha(\Delta (M))^\prime \rightarrow 
 \beta (\Delta (N))^\prime $$ mapping $\cl L_1$ onto $\cl L_2.$
\end{theorem}

\begin{theorem}\label{56}\cite{eleref} Let $M, N$ be CSL algebras (see \cite{dav}). Then $M$ and $N$ are  $\Delta$-equivalent if and only if they are TRO equivalent.
\end{theorem}

\begin{theorem}\label{57}\cite{eleref} Let  $M, N$ be  $\Delta$-equivalent dual operator algebras. For every completely isometric normal representation 
$\alpha $ of $M$ there exists a completely isometric normal representation $\beta $ of $N$ such that the algebras 
$\alpha (M)$ and $\beta(N) $ are TRO equivalent.
\end{theorem}

\begin{theorem}\label{58}\cite{elepure} Two $W\sp{*} $ algebras are 
$\Delta$-equivalent if and only if they are R-WME.
\end{theorem}

\begin{theorem}\label{59}\cite{bk} If two dual operator algebras are 
$\Delta $-equivalent then they are BK-WME.
\end{theorem}

We are going to prove that $\Delta$-equivalence is strictly stronger than BK-WME:

\begin{theorem}\label{510}\cite{elepure, eleref, elenest} Two BK-WME algebras are not necessarily stably isomorphic.
\end{theorem}

\noindent\textbf{Proof} Let $\{q_n: n\in \bb N \}$ be an enumarations of the rationals.
 We define the measure $\mu =\sum_n \delta _{q_n} $ on the Borel 
$\sigma$-algebra of $\bb R,$ where $\delta _{q_n} $ is the Dirac measure. 
Let $H=L^2(\bb R, \mu )$ and let  
$ Q_t^+$ (resp. $Q_t^-$ ) be the projection onto $L^2((-\infty ,t], \mu)$ 
(resp. $L^2((-\infty ,t), \mu)$ ). The set 
$$\cl N_1=\{Q_t^+ , Q_t^-: t\in \bb R \}$$ 
is a totally atomic nest in $H$.  
Let $\lambda $ be Lebesgue measure on $\bb R$ and let $N_t$ be the projection onto $L^2((-\infty ,t], \lambda).$ 
We define the nest 
$$\cl N_2=\{Q_t^+ \oplus N_t,\, Q_t^-\oplus N_t: t\in \bb R \}$$ 
acting on $H\oplus  L^2(\bb R, \lambda  ).$
The map $\theta : \cl N_1\rightarrow \cl N_2$ sending the projection $Q_t^+$ to $Q_t^+\oplus N_t$ and $Q_t^-$ to $Q_t^-\oplus N_t$ 
 is a nest isomorphism \cite{dav}. Therefore, by Theorem \ref{47}, the nest algebras  
$M_1=\mathrm{Alg}(\cl N_1)$ and $M_2=\mathrm{Alg}(\cl N_2)$ 
are BK-WME. Suppose now that they are stably isomorphic. By Theorems \ref{55}
amd  \ref{56}, $M_1$ and $M_2$ should be TRO equivalent. 
So by Theorem \ref{52} there would exist a $*$-isomorphism 
$$\sigma : \cl N_1^{\prime \prime }\rightarrow \cl N_2^{\prime \prime }$$
 such that $\sigma (\cl N_1)=\cl N_2.$ But $\cl N_1^{\prime \prime }$ is a totally atomic maximal abelian selfadjoint algebra  (masa) and  $\cl N_2^{\prime \prime }$ 
is a masa with a nontrivial continuous part. 
This is a contradiction. So $M_1$ and $M_2$ 
are not stably isomorphic. $\qquad \Box$

G.K. ELEFTHERAKIS\\
Department of Mathematics\\Faculty of sciences\\
University of Patras\\265 04 Patras, Greece 

email: gelefth@math.upatras.gr

\end{document}